\documentclass[11pt]{article}
\usepackage{xcolor}
\usepackage{amsmath}
\usepackage{amssymb}
\usepackage{soul} 
\usepackage{pgf}
\usepackage{tikz}
\usepackage{hyperref}
\usepackage{mathrsfs}
\newtheorem{theorem}{Theorem}[section]
\newtheorem{lemma}{Lemma}[theorem]

\newtheorem{remth}{Remark}[theorem]

\newtheorem{conjecture}{Conjecture}

\newcommand{\Frac}{\displaystyle\frac}            
\newcommand{\Int}{\displaystyle\int}              
\newcommand{\Sum}{\displaystyle\sum}              
\newcommand{\refrm}[1]{{\rm(\ref{#1})}}           
\newcommand{\vv}{\vspace{7pt}                     

}
\newcommand{\I}[1]{1\hspace{-3pt} {\rm I}_{#1}}       
\newcommand{\Ib}[1]{1\hspace{-3pt} {\rm I}_{\{#1\}}}  
\newcommand{\Proof}{\noindent {\bf Proof.~}}
\newcommand{\QED}{\hfill \rule{2.5mm}{2.5mm}      

\vspace{7pt}                                        

\noindent }
\newcommand{\df}{:=}
\newcommand{\cC}{\mathscr{C}}

\newcommand{\R}{{\mathbb R}}
\newcommand{\BM}{{\mathbb M}}
\newcommand{\BT}{{\mathbb T}}
\newcommand{\E}{{\mathbb E}}                   
\newcommand{\PR}{{\mathbb P}}                  
\newcommand{\tss}[1]{\textsuperscript{#1}}       
\newcommand{\ca}[1]{{\cal #1}}
\newcommand{\ex}[1]{\,{\rm exp}\set{#1}\,}        
\newcommand{\set}[1]{\left\{#1\right\} }          
\newcommand{\ov}[1]{\overline{#1}}
\newcommand{\imp}{\,\Rightarrow\,}

\newcommand{\vspandexsmall}{\rule[-11pt]{0pt}{22pt}}   

\DeclareMathOperator {\Var}{Var}
\DeclareMathOperator{\im}{imp}
\newcommand{\equ}{\,\Leftrightarrow\,}
\DeclareMathOperator{\sinc}{sinc}
\DeclareMathOperator{\sign}{sign}
\DeclareMathOperator{\Ei}{Ei}

\DeclareMathOperator{\eve}{eve}

\newcommand{\be}{\begin{equation}}
\newcommand{\ee}{\end{equation}}

\newcommand{\footer}[1]{{\def\thefootnote{}\footnotetext{#1}}}
\begin{document}\sf
\thispagestyle{empty}

\centerline{\bf Improper Poisson integral of the sinc function}
\vv
\begin{center}
Jerzy Szulga{$^\dagger$}\footer
{\begin{tabular}{l}
$^\dagger$ Department of Mathematics and Statistics, Auburn University, USA\\
{~} MSC 2020: Primary 60H06. Secondary 60F25\\
{~} Key words and phrases: improper Poisson integral, renewal process, LePage series, \\
\hspace{104pt} sinc function
\end{tabular}}
\vv
\begin{minipage}{0.9\linewidth}
{\small {\bf Abstract.} A  formal sum  $\sum_n f(S_n)$ may be seen as the integral $\int f dN$ with respect to random point process $N(A)=|\{n:S_n\in A\}|$. We study its convergence beyond the well known context of Lebesgue integrable functions, admitting nonintegrable functions whose improper Riemann-Lebesgue integrals exist. We focus on the sinc function and some of its relations leaving the general case to conjectures.}
\end{minipage}
\end{center}

{\footnotesize \tableofcontents}

\section{Rationale}
A sequence of random elements $S_n$ in a measure space $(\BM,\ca M,\mu)$ entails the random counting measure $N(A)=|\{n: S_n\in A\}|$ on $\ca M$ which in turn induces a linear operator on some class $\cC$ of measurable functions 
\be\label {Nf}
N f=\sum_n f(S_n),
\ee
where the series converges at least in probability if not a.s. or in a stronger metric.
It is well defined for nonnegative function and for simple functions $f$  we recognize the integral $Nf=\int_{\BM} f\,dN$   although beyond $\BM=\R$ the paths ``$N_t$'' are not well defined. Not much is known about further extensions  beyond the well known case of a Poisson random measure and (almost) integrable functions. That is, the adverb ``almost'' refers to the condition $\Psi_1(f)\df \Int_{\BM} (|f|\wedge 1)\,d\mu<\infty$ that defines a modular space $L^{\Psi_1}$. The condition is sufficient for the a.s.\ convergence of the series in \refrm{Nf}  (cf., e.g., \cite{KalS,Mus}).
Indeed, this can be deduced from the Laplace transform 
\be\label{L}
\E e^{-N f} =e^{-\mu\left(1-e^{-f}\right)},\quad f\ge 0,
\ee
Here and in the sequel we use the operator notation for an integral, e.g., $\mu f=\int_{\BM}f\,d\mu$ or with the Lebesgue measure $\lambda f=\int_0^\infty f(x)\,dx$, etc. 
\vv

In particular, the a.s. convergence of series \refrm{Nf} necessarily requires that $\mu\{|f|>c\}<\infty$ 
so $N[ f;|f|>c]$ is but a finite sum.
Therefore, in a further study w.l.o.g.\ we may and do confine to bounded functions, and thus the integrability becomes the sufficient condition for the a.s.\ convergence of series \refrm{Nf}.
Yet, for a Poisson measure the ch.f.\ is well defined assuming only the existence of the Riemann-Lebesgue improper integral:
\be\label{F}
\E e^{\imath \, Nf} =\ex{-\mu \left(1-e^{\imath f}\right)}
=\ex{- \Big(  \mu \left(1-\cos  f \right)  +  \imath\,\mu \sin  f \Big)}.
\ee
Hence, by symmetrization, $(N-N')f$ exists and the series in  $Nf^2$ converges a.s. Beyond a Poisson measure such inference my be incorrect. The improper integral is understood as usual: for some sequence  $A_n\in \ca M$ such that $\mu A_n<\infty$, $A_n\nearrow \BM$, and 
\[
\mu\, \sin f \df \lim_n \int_{A_n} \sin f\,d\mu,
\]
Thus we arrive at a Banach space with an ``impaired'' integrability: 
\be\label{Limp}
L^2_{\im}= L^2\cap L^\infty\cap \{\text{improper integrals $\mu f$ exist}\}
\ee
with a norm, e.g., $\|f\|_2+\|f\|_\infty+|\mu f|$. Since $L^1\cap L^\infty$ (or the subspace of simple functions with supports of finite measure) is dense in $L^2_{\im}$, the standard completeness argument entails the well defined \cite{KwaW} Wiener-type stochastic integral $\int_{\BM} f\,dN$  on $L^2_{\im}$.   
However, its mere existence does not immediately guarantee that the series on the right hand side of \refrm{Nf} is a ``real thing'', i.e., that its partial sums converge even in probability. Further, it is not presently clear whether the terms $f(S_n)$ are summable (i.e., the series converges unconditionally). 
\vv
Functions $f$ and $g$ are dubbed equivalent if $f-g\in L^{\Psi_1}$ eventually; so we write $f\smile g$.
\begin{conjecture} \label{conj:as} Consider a Poisson random measure.
\begin{enumerate}
\item On $L_{\im}$ series \refrm{Nf}  converges a.s.\ and in $L^2$ eventually  (i.e., for $n\ge n_0$).
\item The convergence is not unconditional.
\end{enumerate}
\end{conjecture}
While we have not been able to prove or disprove it in general, nevertheless we investigate the sinc function in detail,  $f(x)=\sinc x=\Frac {\sin x} x$, a flag member of $L_{\im}$, and some of its relations.

\section{The sinc function}
\subsection{A renewal process}
Let $(X,X_k)$ be a sequence of i.i.d.\  nonnegative random variables. Assuming also a nondegenerate distribution, denote $z=\E e^{\imath X}=re^{\imath \alpha}$, $r<1$. For parametrized quantities $A$ and $B$ we write $A\prec B$ if for some constat $c>0$, independent of the parameters, $A\le cB$; also  $A\asymp B$ if $A\prec B$ and $B\prec A$. Denote 
\[
R_n=\frac 1 {S_n},\quad \zeta_n=e^{\imath S_n},
\quad D_n=R_n-R_{n+1},\quad Z_n=\sum_{k=1}^n \zeta_k,\,Z_0=0.
\]
\begin{theorem}\label{th:If} Let $f(x)=\Frac{e^{\imath x}}{x}$. 
Then series \refrm{Nf} converges a.s. If $\E X^{2p}<\infty$ then the series converges eventually in $L^p$ 
for every $p> 0$. The convergence is unconditional for each of these modes.
\end{theorem}
We need several technical details prior to the proof. 

\begin{lemma} \label{mart} 
Let $X_k$ be independent copies of a nondegenerated random variable $X$.
\begin{enumerate}
\item $\|Z_n\|_2\asymp \sqrt{n}$. 
\item 
$Z_N=\Sum_{k=0}^N z^{N-k} M_k$ for some  martingale $M_k$.
\item $\|Z_n\|_p\asymp \|Z\|_2\asymp \sqrt{n}$ for every $p\ge 2$ (with constants dependent on $p$).
\end{enumerate}
\end{lemma} 
\Proof First we estimate the $L^2$-norm. We have
\[
z_{k}\df \E Z_{k} =\sum_{\ell=1}^{k} z^\ell=\frac{z}{1-z}(1-z^{k})=\Frac z {1-z} -\Frac{z^{k+1}}{1-z}, \,\,k\ge 1,
\]
and $z_0=0$. Then, for $k\ge 1$
\[
\E|Z_k|^2=\E|1+Z_{k-1}|^2=1+z_{k-1}+\ov{z}_{k-1} +\E|Z_{k-1}|^2.
\]
Summing up and putting $c=c_z=1+2\Re \Frac z{1-z}$ we obtain $\E|Z_n|^2=c\,n +O(1)\asymp n$.
\vv
Next, we set up Doob's decomposition for the standard filtration $\ca F_n$:
\[
Z_n=M_n+A_n, \quad M_n=\sum_{n=1}^n \Big(\Delta Z_n -
\E\left[\Delta Z_n|\ca F_{n-1}\right]\Big),\quad
A_n=\sum_{n=1}^n  \E\left[\Delta Z_n|\ca F_{n-1}\right].
\]
Put $Z_0=M_0=A_0=0$.
Then by straightforward computation we check that for $n\ge 1$
\[
\E\left[\Delta Z_n|\ca F_{N-1}\right] = z\,\Delta Z_{N-1}
\quad \imp \quad Z_n=M_n+z\,Z_{N-1},
\]
and the claimed convolution follows by iteration.
\vv
By the BDG inequalities \cite{BDG}, $\|M_n\|_p\asymp \|Q_n\|_p$, where $Q_n^2=\sum_{k=1} |\Delta M_k|^2\asymp n$ (with 1 as the upper constant and $\frac {1-r^2}{1+r^2}$ as the lower constant).
This ensures the equivalence for $p\ge 2$. Indeed,  
\[
\|Z_n\|_p \le \sum_{k=1}^n r^{n-k} \|M_k\|_p \le C_p \sum_{k=1}^n r^{n-k} \sqrt{k}\asymp \sqrt{n}\asymp \|Z_n\|_2.
\]
since $\sum_{k=1}^n r^{-k}\sqrt{k}\asymp r^{-n}\sqrt{n}$. 
\QED
\begin{conjecture}
$\|Z_n\|_p\asymp \|Z_n\|_2$ for $p<2$. 
\end{conjecture}
Yet, we will not need this relation in the sequel; the obvious dominance ``$\prec$'' suffices.

\begin{lemma} \label{LpR}
$\|R_n\|_p\prec n^{-1}$ eventually.
\end{lemma}
\Proof  Although for small $n$ the left hand side might be infinite yet $\E X^{-\alpha}<\infty$ for $\alpha<1/2$. Indeed, 
w.l.o.g.\ (with a little proof\footnote{e.g., replace $X:=X+\lambda V$, where $V$ is exponential and independent of $X$, and with the obtained density complete the argument and let $\lambda\to 0$, utilizing the Lebesgue's monotone convergence }) we may and do assume that $X$ admits a density $g(x)$ and then use the  CSI.
So, if $2\alpha<1$,
\[
\E [X^{-\alpha};X\le 1] =\int_0^1 \frac{g(x)}{x^\alpha}\le \left(\int_0^1 g^2(x)\,dx\right)^{1/2}
\left(\int_0^1 x^{-2\alpha}\,dx\right)^{1/2}\le (1-2\alpha)^{-1}
\]
Hence, comparing the arithmetic and geometric means and denoting $X_n^*=\max_k X_k$,
\[
\E R_n^p \Ib{X^*\le 1}\le \frac 1 {n^p}\, \E \frac {\Ib{X_n^*\le 1}}{(X_1\cdots X_n)^{p/n}}
=\frac {\left(\E[X^{-p/n}; X^*_n\le 1]\right)^n}{n^p} <\infty
\]
when $n>2p$. So, by the Lebesgue's dominated convergence and then by the SLLN,
\[
\limsup_n n^p \, \E R_n^p\le \E \limsup_n \left(\Frac{n}{S_n}\right)^p=\Frac1 {(\E X)^{p}},
\]
which makes sense even when $\E X=\infty$.
\QED
{\bf Proof of Theorem \ref{th:If}. } 

Using Abel's change of summation:
\be\label{Abel}
\Sum_{n=1}^N R_n D_n =R_NZ_N-\Sum_{n=1}^N D_n Z_n\\
\ee
First, by algebra and \ref{LpR},
\[
D_n \le X_{n+1} R_n^2\quad\imp\quad \E D_n^2\le \frac{E X^2}{n^4}
\]
Now, we will prove the a.s.\ convergence. The series $\Sum_{n=1}^N D_n Z_n$ converges absolutely in $L^1$ which follows from  the CSI and \refrm{LpR}:
\[
\E|D_nZ_n|\le \|D_n\|_2\|Z_n\|_2 \prec \frac 1 {n^{3/2}}.
\]
Hence the series converges a.s. In the first term we separate $R_N Z_N=\Frac N {S_N} \cdot \Frac{Z_N}{N}$. The first ratio is subject to the standard SLLN. Although the summands $\zeta_n$ in the second ratio lack even a trace of independence, they are bounded and their sum has a relatively small second moment. By \cite[Thm.1]{Lyo} the condition
\[
\sum_N \Frac{\|Z_N\|_2^2}{N^3}<\infty
\]
is sufficient for the SLLN with 0 limit, and it is satisfied since $\|Z_n\|_2^2\asymp n$. 
\vv
Now, the first single term in \refrm{Abel} satisfies 
\[
\|R_nZ_n\|_p \le  \|R_n\|_{2p}  \,   \|Z_N\|_{2p}\prec \Frac 1 {\sqrt n}.
\] 
Again by CSI, in virtue of Lemma \ref{LpR},
\[
\|D_nZ_n\|_p\le \|X\|_{2p}\,\|R_n\|_{4p}^2\,\|Z_n\|_{2p}\prec \frac{1}{n^{3/2} }, 
\]
so the sum in \refrm{Abel} eventually converges  absolutely in $L^p$. 
\vv
Finally, we observe that the applied tools (absolute convergence of the series stemming from \refrm{Abel} and the sufficient boundedness of the summands in  \cite[Thm.1]{Lyo}) work equally well for permutations of the original sequence $S_k$.
\QED
\subsection{Alternative atomic representations}
\subsubsection{LePage series} 
Although arrival moments $S_n$ seem to be natural atoms of a Poisson process on $\R_+$ but  atomic representations are not unique as illustrated by Poisson measures tantamount to positive pure jump L\'evy processes $X_t$ on $\R_+$.
That is, for a one-sided L\'evy measure $\nu$,
\be\label{LT}
\E e^{-\theta X_t} =e^{-t \psi(\theta)}. \quad\text{where}\quad \psi(\theta)
=\int_0^\infty \left(1-e^{-\theta x}\right)\,\nu(dx)
\ee
A finite integral exists iff $\nu(x\wedge 1)<\infty$. 
When $\R_+$ is replaced by an arbitrary Borel space $(\BT,\ca T,\tau)$, the paths ``$X_t$'' may be undefined. 
Yet, the L\'evy integral $X f\df\int_{\R_+} f\,dX$ for a Borel function $f\ge 0$ on $\BT$  has the Laplace transform
\be\label{XfT}
\ln \E e^{-X f} =(\tau\otimes \nu) \psi(f)
=\int_{\BT} \int_0^\infty \left(1-e^{-f(t) x}\right)\,\tau(dt)\, \nu(dx),
\ee
We note the correspondence to the abstract Poisson measure \refrm{L} on $\BM=\BT\times \R_+$ with the control measure $\mu=\tau\otimes\nu$, and the function $F(x,t)=xf(t)$. By the Borel isomorphism we can reduce $\BT$ to $\R_+$ and $\tau$ to the Lebesgue measure $\lambda$, so we will assume the latter henceforth.
\vv

We consider now LePage-style series representations (see \cite{LeP}, utilized there for stable processes, expanded to the present form in \cite{Ros}).  Denote $G(x)=\nu(x,\infty)$ and $H=G^{-1}$. 
Consider independent copies  $V_n$  of a random variable $V$ with a strictly positive density $p(x)$  on $\R_+$ that are independent of the Poisson arrivals $S_n$.  put $H_n=H(S_n\,p(V_n)$. First, define the ``atomic'' version akin to \refrm{Nf}:
\be\label{H}
X f\df \sum_n H_n \,f(V_n),
\ee
initially for $f\ge 0$.  
Referring to the discussion around formula \refrm{Nf}, the series converges a.s.\ on the modular space \cite{Mus}:
\be\label{Lpsi}
L^{\Psi}:\quad\left\{\begin{array}{rcl}
\Psi_2(f) &=& \Int_{\R_+} \Int_0^\infty \left(x^2f^2(t)\wedge 1 \right)\,\nu(dx)\,dt<\infty,\\
\Psi_1(f) &=& \Int_{\R_+} \Int_0^\infty (|xf(t)|\wedge 1)\,\nu(dx)\,dt<\infty
\end{array}\right.
\ee
That is, the stochastic processes $Xf$ and $\int_{\R_+} f\,dX$ have the same distribution over $L^\Psi$. 
\vv

The general ch.f.\ \refrm{F}, adjusted for the product measure as in \refrm{XfT} (with $\BT=\R_+$), becomes  
\be\label{Xfpsi}
\E e^{\imath\,X f} =\ex{\lambda \psi_r(f)+\imath\,\lambda \psi_i(f)},
\ee
involving two functions on $\R$:
\be\label{2int}
\begin{array}{rcl}
\vspandexsmall
\psi_i(\theta) &=& \Int_0^\infty \sin(\theta x)\,\nu(dx),\\
\psi_r(\theta) &=& \Int_0^\infty \big(1-\cos (\theta x)\big)\,\nu(dx).
\end{array}
\ee
Since the former function directs us  beyond the modular space \refrm{Lpsi} thus the question of improper integral emerges again. 
\vv

Conditioning upon the Poisson arrivals $(S_n)$ we deal with a sum of independent random variables whose a.s.\ convergence is controlled by the Kolmogorov's Three Series Theorem, and further, by three Poisson integrals, proper or improper. Namely, denoting $[x]_c\df x\Ib{|x|\le c}$:
\be\label{KKK}
\begin{array}{rcl}
K_1(s) &=& \PR(H(sp(V)) |f(V)| >1),\\
K(s)=K_2(s) &=& \E[H(sp(V)) f(V)]_1,\\
K_3(s) &=& \E[H^2(sp(V)) f^2(V)]_1,\\
\end{array}
\ee
with any constant $c>0$ replacing the cut-off 1
and $K_3-K_2^2$ yielding the variance. Thus, the series in \refrm{H} converges a.s.\ iff the three Poisson integrals $N K_i$, $i=1,2,3$, converge a.s. Only $N K$ may be improper. Since $K_1$ and $K_3$ are bounded, the related convergence is independent of the choice of the density $p(v)$ and reduces to the integrability of $K_1$ and $K_3$:
\[
\begin{array}{rl}
N K_1: & \Int_0^\infty G(1/|f(v)|) \,dv=\Int_0^\infty \Int_{1/|f(v)|}^\infty\,\mu(du) \,dv <\infty \\
N K_3: &
\Int_0^\infty f^2(v) \Int_{1/|f(v)|} ^\infty u^2\,\nu(du)\,dv<\infty.
\end{array}
\]
The integrals are finite iff the modular $\Psi_2(f)$ is finite, which is guaranteed by the existence of $X f^2$. Therefore, the potential issue of improper integration has been reduced to the case of a Poisson process on $\R_+$ and the function $K$. 

\vv
\subsubsection{Examples}
Although the new function $K(s)$ is more complicated than the original function $f$ yet it may improve its properties for selected densities $p$. Only large $s$ are relevant, so in further considerations we assume $s\gg1$. Hence we employ the mark ``eve''(for ``eventually'') in writing $L^1_{\eve}$ or $L^2_{\eve}$. 
\vv

For a standard Poisson process $N_t$ with Poisson arrivals $S_n$
\[
\psi_r(\theta)=1-\cos\theta,\quad \psi_i(\theta)=\sin\theta.
\]
Notice that  $\psi_i(f)\smile f$ on $L_{\im}$.
Since $\nu=\delta_1$ is the L\'evy measure of the standard Poisson process, hence $G=H=\I{[0,1]}$. In this case, $\lambda K=\lambda f$ exist or do not exist simultaneously.  
A choice of a monotonically decreasing density entails the inverse $q=p^{-1}$. Further, we may assume that a chosen function $p>0$ that is just integrable. Thus
\be\label{plainK}
K(s)=\int_a f(v) p(v)\,dv\quad\text{with}\quad a=a(s)={q(s^{-1})}.
\ee
Because of a  computational feasibility we consider $p(x)=x^{-r}, \,r>1, \,x\gg0$, stemming from a Pareto density. That is,
\[
K(s)=\int_{a} ^\infty f(v)x^{-r} \,dv ,\quad a =s^{1/r}.
\]
For the sinc function, immediately integrating by parts,
\[
K(s)=\int_{a} ^\infty \Frac{\sin v}{v^{r+1}} \,dv=\Frac{\cos a}{a^{r+1}}-\int_a^\infty \Frac{\cos v}{v^{r+2}}\,dv
\]
i.e., $K(s)\asymp s^{-1-1/r}$, a function eventually integrable.
\vv

A similar argument applies to $f(x)=\frac{\cos x} x$ eventually, or even for slightly more general functions $f$.
\begin{theorem}
For a Poisson process, there is a density $p(x)$ such that series \refrm{H} converges a.s.\ for a complex valued function $f(x)=A(x)\sin x$ with an amplitude $A\in L^2_{\eve}$ such that  $A'\in L^1_{\eve}$. 
\end{theorem}
\Proof We choose a Pareto density like above.
That is, with $a=s^{1/r}$ our function $K(s)$ becomes
\[
\begin{array}{rcl}
\Int_a^\infty \Frac{A(x)e^{\imath x}} {x^r}\,dx &=& \imath \left( \Frac{A(a) e^{\imath a}}{a^{r}}+
\Int_a^\infty \Frac{A'(x)e^{\imath x}} {x^r}\,dx -
 r \Int_a^\infty \Frac{A(x)e^{\imath x}} {x^{r+1}}\,dx
 \right)\\
\end{array}
\]
Therefore we obtain $K \in L^1_{\eve}$:
\[
\begin{array}{rcl}
\text{1\tss{st}  quotient}: &&\text{by CSI since  $A\in L^2_{\eve}$ and } \,\Frac {A(s^{1/r})}{s}\in L^1_{\eve}\equ 
\Frac {A(x)}{x}\in L^1_{\eve};\\
\text{1\tss{st} integral}: && A'\in L^1_{\eve};\\
\text{2\tss{nd} integral}: &&\text{by CSI since $A\in L^2_{\eve}$},\\
\end{array}
\]
which completes the proof.
\QED
\begin{remth}\sf
For the exponential density $p(x)=e^{-x}$ and  the sinc function, \refrm{plainK} yields
\[
K(s)=\int_{a} ^\infty \sin v \,\frac{e^{-v}\,dv}{v},\quad a(s)=\ln s.
\]
Computing (cf., e.g., the Wolfram integral calculator \cite{Wolf}),
\[
K(s)=\iota(a)=\frac \imath 2 \Big(Ei(-1+\imath)a) - \Ei ((-1-i)a)  \Big)
\smile e^{-a}\sinc (a)=\frac{\sin(\ln s)}{s\ln s},
\]
and so we are back in square one. The integration by parts that worked for a Pareto density moves us around in circles. Yet, the ``integrability'' slightly improved, not enough though, and the $L^2$-norm of $K(s)$  is much smaller than the $L^2$-norm  of $\sinc s$.
\end{remth}
\section{Discussion}
\subsection{Other L\'evy processes}

For a positive $\alpha$-stable process $X_t$, $\alpha<1$ the L\'evy measure
$\nu(x)=\Frac \alpha {x^{\alpha+1}}\, dx$, so
\[
\begin{array}{rcll}
\psi_r(\theta) &=&c_\alpha |\theta|^p,              & c_\alpha =\Gamma(1+\alpha)\cos(\alpha\pi/2),\\
\psi_i(\theta) &=& s_\alpha \sign(\theta)|\theta|^\alpha,& s_\alpha =\Gamma(1-\alpha)\sin(\alpha\pi/2)
\end{array}
\]
(cf.\ \cite[3.761.4]{GR} with standard calculations). Here $Xf$ exists iff $f\in L^\alpha$ and the issue of  improper integration does not even exist.
\vv

For a Gamma process $X_t$ (a fractional generalization of a Poisson $N_t$):
\[
\psi_r(t)=\frac 1 2\, \ln(1+\theta^2),\quad \psi_i(t)=\arctan\theta
\]
(cf., e.g., \cite{Szu:gam}). Then $G(x)=E_1(x)=\Int_x \Frac{e^{-u}} u\,du $, the exponential-integral function.
On $L^2_{\im}$, since $\psi_i( f)\smile f$  the issue of the improper integration is valid. 
However, we haven't been able to establish the integrability of the function $K(s)$ \refrm{KKK} even for previously successful densities $p(x)$ and even for the specific sinc function.
\vv

We note that the class of such functions is scale and shift invariant as $L^2_{\im}$ is. However, this property is not immediately satisfied by $\cC$ (defined above \refrm{Nf})). 
 So, we formulate also a weaker
\begin{conjecture}
The space $\cC$ is invariant under affine transformations of the variable.
\end{conjecture}

\subsection{A potentially misleading simulation}
It may be tempting to rely on synthetic series \refrm{Nf} that are easy to simulate. However, a simulation must be necessarily restricted to finite sums which may be a source of confusion because it produces ``visually'' convergent ``series''.
\vv
The following approach gives some sense to the series but doesn't answer the question, implicit in Conjecture \ref{conj:as}. That is, partitioning the halfline into disjoint intervals $I_k$  yields  $f_k=f\I{I_k}$ and  independent random variables $X_k=N f_k=\sum_n f_k(S_n)$ a.s. Then the series
\[
\sum_k \E X_k=\sum_k \lambda f_k=\lambda f \quad \text{and}\quad \sum_k \Var(X_k)=\sum_k \lambda f_k^2 \quad\text{converge},
\]
hence by Kolmogorov's Two Series Theorem,
\be\label{eq:strong flawed}
\text{``$N f$''}\df\sum_k X_k =\sum_k \left(\sum_n  f_k(S_n)\right) \quad\text{converges a.s.\ and in $L^2$}. 
\ee
A simulation of series \refrm{Nf} must use a finite $n$-sum in which case the above order of summation can be reversed but not so for the infinite $n$-sum when $f$ is not integrable, at least not obviously.  Thus, we exit with yet another conjecture;  recalling \refrm{Limp}.
\begin{conjecture}
Let $f\in L^2_{\im}$ and $A_k$ be disjoint Borel subset of $\R_+$. Let $S_n$ be Poisson arrivals.  Then the double sequences $f_k(S_n)$ is  summable a.s.
\end{conjecture}

\addcontentsline{toc}{section}{References}
\vv

{\footnotesize  
\begin{thebibliography}{9}

\bibitem{BDG}
D.L. Burkholder, B. Davis, R.F. Gundy, 
Integral inequalities for convex functions of operators on martingales, 
Proc. 6th Berkeley Symp. Math. Statistics and Probability, 2 (1972) pp. 223–240.

\bibitem{GR}
I.S.\ Gradshteyn and M.\ Ryzhik,  {\it Table of Integrals, Series, and Products}. Academic Press, Inc., San Diego, CA, 1979.

\bibitem{KalS}
O.\ Kallenberg and J.\ Szulga, Multiple integration with respect to Poisson and L\'evy processes, Prob.\ Th.\ Rel.\ fields, 83, 101-134, (1989).

\bibitem{KwaW}
S.\ Kwapie\'n and W.A.\ Woyczy\'nski, {\em Random Series and Stochastic Integrals: Single and Multiple}, Birkh\"auser, Boston 1992 (2nd ed. 2000).

\bibitem{LeP}
R.\ LePage, Multidimensional infinitely divisible variables and processes. Part I: Stable case, Technical Report, Statistics Department, Stanford (reprinted in Lect.\ Notes Math., Vol. 1391 (1989), p. 292).

\bibitem{Lyo}
R.\ Lyons,
Strong law of large numbers for weakly correlated random variables.
 Michigan Math. J. Vol.35, No. 3, 352-359 (1988).

\bibitem{Mus}
J.\ Musielak, {\em Orlicz Spaces and Modular Spaces}, Lect. Notes Math. 1034, Springer, Berlin 1986.

\bibitem{Ros}
J.\ Rosi\'nski, On series representation of infinitely divisible random vectors. Ann. Probab. 18 (1990), 405–430.

\bibitem{Szu}
J.\ Szulga, {\em Introduction to Random Chaos}, Chapman \& Hall, London 1998.

\bibitem{Szu:gam}
J.\ Szulga, Random Gamma time,  Asian J.\ Stats.\ Appl.), 1 (1) (2024), 31-57. 


 


\bibitem{Wolf} 
WolframAlpha Online Integral Calculator.

\url{https://www.wolframalpha.com/calculators/integral-calculator}

\end{thebibliography}}
\end{document}